# Multi-objective Optimization: A Case Study

Nazmul Hasan

## Abstract

The aim of this literature is to illustrate the application of multi-objective optimization routines through a case study of face milling operation. For this purpose, the face milling operation is designed as a multi-objective optimization problem and then solved to obtain optimum values for the machining parameters– cutting speed ($Vc$), feed rate ($fz$) and depth of cut ($t$); using the optimization routines. The formulated problem of face milling operation includes two conflicting objectives– to maximize Material Removal Rate ($MRR$) and to minimize surface roughness ($Ra$). Among various multi-objective optimization routines, five of them namely Global Criterion Method, Lexicographic Method, Weighted Sum Method, ε – Constraint Method and Genetic Algorithm are used in this literature. The outcomes of these multi-objective optimization routines are then compared to reflect their relative attractiveness.

## Nomenclature

| | | | | | |
|---|---|---|---|---|---|
| $Vc$ | : | Cutting speed, m/min | $h_k(x)$ | : | $k$th equality constraints |
| $fz$ | : | Feed rate, mm/tooth | $lb$ | : | Lower bound |
| $t$ | : | Depth of cut, mm | $ub$ | : | Upper bound |
| $Ra$ | : | Surface Roughness, μm | $p$ | : | Parameter in Global Criterion Method |
| $MRR$ | : | Material Removed Rate, mm³/min | $w_i$ | : | Weighing factor |
| $RSM$ | : | Response Surface Methodology | $\varepsilon_i$ | : | Upper bound for $i$th objective function |
| $APD$ | : | Absolute Percentage Deviation | $\varepsilon_{Ra}$ | : | Upper bound for Surface Roughness, μm |
| $MAPD$ | : | Mean Absolute Percentage Deviation | $f_i^{norm}$ | : | Normalized Value of $i$th objective function |
| $x$ | : | Design variables | $f_i^{min}$ | : | Minimum value of $i$th objective function |
| $f_i(x)$ | : | $i$th Objective function | $f_i^{max}$ | : | Maximum value of $i$th objective function |
| $f_i(x^*)$ | : | Optimum value of $i$th objective function | $Ra_{min}$ | : | Minimum value of Surface Roughness, μm |
| $g_j(x)$ | : | $j$th inequality constraints | $Ra_{max}$ | : | Maximum value of Surface Roughness, μm |

# 1.0 Introduction

In many practical applications, designer may want to optimize two or more objective functions simultaneously. These are called multi-objective, multi-criteria or vector optimization problems (Arora, 2012). Multi-objective optimization techniques can play an important role in many areas of machining processes, both conventional and modern (Rao et al., 2017). Therefore, a case study of face milling has been adopted in this literature from the paper titled "Modeling and optimization of surface roughness and productivity through RSM in face milling of AISI 1040 steel using coated carbide inserts" by Fnides et al., 2017. In this literature, models developed using RSM have been proposed for the response variables of interest and optimization has been performed in the paper using desirability approach.

In the current literature, different multiple non-linear regression models have been proposed for the response variables of interest and suitability of the models are compared with those have been mentioned in the paper. For optimization, other than desirability approach, several multi-objective optimization techniques have been used, which are−

  i. Global Criterion Method
 ii. Lexicographic Method
iii. Weighted Sum Method
 iv. $\varepsilon$ – Constraint Method
  v. Genetic Algorithm

Among these methods, Weighted Sum Method, $\varepsilon$ – Constraint Method and Genetic Algorithm can be used to produce Pareto optimal set and generate Pareto optimal front (Messac et al. 2000). Pareto optimal front can help the decision makers to choose a suitable tradeoff between the conflicting objectives in a multi-objective optimization problem.

# 2.0 Optimum Design Formulation:

Proper formulation of a design optimization problem is important because the optimum solution is as good as the formulation. Usual practice for most design optimization problems is five-step formulation procedure (Arora, 2012). The five steps are:

Step 1: Project/problem description

Step 2: Data and information collection

Step 3: Definition of design variables

Step 4: Optimization criterion

Step 5: Formulation of constraints

In this literature, these five steps will be followed to address the multi-objective optimization in the case study of face milling operation.

## 2.1 Project/problem description

In metal machining process, like face milling operation, high Material Removal Rate (*MRR*) is desired to achieve high productivity. High MRR involves high values of machining parameters – cutting speed (*Vc*), feed rate (*fz*) and depth of cut (*t*). On contrary, high values of machining parameters cause high surface roughness (*Ra*), which is undesirable. Therefore, *MRR* is to be maximized while *Ra* needs to be minimized. Machining parameters need to be selected to ensure suitable tradeoff among the conflicting machining responses−Material Removal Rate (*MRR*) and surface roughness (*Ra*). The machining parameters – cutting speed (*Vc*), feed rate (*fz*) and depth of cut (*t*) must be within the following limits:

$$78 \leq \text{cutting speed } (Vc) \leq 314 \text{ m/min}$$
$$0.04 \leq \text{feed rate } (fz) \leq 0.16 \text{ mm/tooth}$$
$$0.2 \leq \text{depth of cut } (t) \leq 0.6 \text{ mm}$$

## 2.2 Data and information collection

For face milling operation, in order to establish functional relationship between the machining responses and machining parameters for face milling, data need to be collected by conducting experiments. Experimental data those has been collected from the paper titled "Modeling and optimization of surface roughness and productivity through RSM in face milling of AISI 1040 steel using coated carbide inserts" (Fnides et al., 2017) are shown in Table 2.1. The paper under study includes the mathematical models expressed by Equations (2.1) and (2.2) for *Ra* and *MRR* respectively.

$$
\begin{aligned}
Ra = &\ 2.65599 - 0.00726733 \times Vc + 1.70439 \times fz - 0.012765 \times t - 0.00273646 \times Vc \times fz \\
&+ 0.000505119 \times Vc \times t + 1.47321 \times fz \times t
\end{aligned} \quad (2.1)
$$

$$
\begin{aligned}
MRR = &\ 7927.21 - 43.31 \times Vc - 84934.4 \times fz - 19818 \times t + 464.12 \times Vc \times fz \\
&+ 108.295 \times Vc \times t + 212917 \times fz \times t
\end{aligned} \quad (2.2)
$$

For further investigation under the project, different multiple non-linear regression models have been proposed. The new mathematical models for *Ra* and *MRR* are expressed in Equations (2.3) and (2.4) respectively.

$$
\begin{aligned}
Ra = &\ 3.082776 - 0.01425 \times Vc + 4.330794 \times fz + 0.465279 \times t + 1.85 \times 10^{-5} \times Vc^2 \\
&- 5.78704 \times fz^2 - 0.15278 \times t^2 - 0.00862 \times Vc \times fz - 0.00142 \times Vc \times t - 2.73511 \times fz \times t \\
&+ 0.018687 \times Vc \times fz \times t
\end{aligned} \quad (2.3)
$$

$$\begin{aligned} MRR = \quad & -2345.09 + 10.37705 \times Vc + 24983.71 \times fz + 7079.912 \times t - 0.01672 \times Vc^2 \\ & - 34722.2 \times fz^2 - 4166.67 \times t^2 - 64.9643 \times Vc \times fz - 13.6425 \times Vc \times t - 66322.5 \times fz \times t \\ & + 1403.922 \times Vc \times fz \times t \end{aligned} \qquad (2.4)$$

Table 2.2 shows experimentally obtained values of *Ra* and *MRR*, the output of the two sets of models ([Eqs. 2.1 & 2.2] and [Eqs. 2.3 & 2.4]) and calculated Absolute Percentage Deviation (APD) values for both sets of models.

**Table 2.1: Experimental data for face milling operation of AISI 1040 steel**

| Experimental Runs | Machining parameters | | | Machining responses | |
|---|---|---|---|---|---|
| | $Vc$ (m/min) | $fz$ (mm/tooth) | $t$ (mm) | $Ra$ (μm) | $MRR$ (mm$^3$/min) |
| 1 | 78 | 0.04 | 0.2 | 2.23 | 730 |
| 2 | 78 | 0.04 | 0.4 | 2.29 | 1460 |
| 3 | 78 | 0.04 | 0.6 | 2.32 | 2190 |
| 4 | 78 | 0.08 | 0.2 | 2.37 | 1460 |
| 5 | 78 | 0.08 | 0.4 | 2.4 | 2920 |
| 6 | 78 | 0.08 | 0.6 | 2.42 | 4380 |
| 7 | 78 | 0.16 | 0.2 | 2.58 | 2920 |
| 8 | 78 | 0.16 | 0.4 | 2.6 | 5840 |
| 9 | 78 | 0.16 | 0.6 | 2.62 | 5760 |
| 10 | 157 | 0.04 | 0.2 | 1.5 | 1460 |
| 11 | 157 | 0.04 | 0.4 | 1.54 | 2920 |
| 12 | 157 | 0.04 | 0.6 | 1.55 | 4380 |
| 13 | 157 | 0.08 | 0.2 | 1.59 | 2920 |
| 14 | 157 | 0.08 | 0.4 | 1.6 | 5840 |
| 15 | 157 | 0.08 | 0.6 | 1.61 | 8760 |
| 16 | 157 | 0.16 | 0.2 | 1.62 | 5840 |
| 17 | 157 | 0.16 | 0.4 | 1.64 | 11680 |
| 18 | 157 | 0.16 | 0.6 | 1.65 | 17520 |
| 19 | 314 | 0.04 | 0.2 | 0.5 | 2920 |
| 20 | 314 | 0.04 | 0.4 | 0.48 | 5840 |
| 21 | 314 | 0.04 | 0.6 | 0.51 | 8760 |
| 22 | 314 | 0.08 | 0.2 | 0.55 | 5840 |
| 23 | 314 | 0.08 | 0.4 | 0.6 | 11680 |
| 24 | 314 | 0.08 | 0.6 | 0.62 | 17520 |
| 25 | 314 | 0.16 | 0.2 | 0.65 | 11680 |
| 26 | 314 | 0.16 | 0.4 | 0.76 | 23360 |
| 27 | 314 | 0.16 | 0.6 | 0.82 | 35040 |

**Table 2.2: Comparison of Models' Results**

| Experimental Values | | Models of Eqs. (2.3) and (2.4) | | Models of Eqs. (2.1) and (2.2) | | APD for Models of Eqs. (2.3) and (2.4) | | APD for Models of Eqs. (2.1) and (2.2) | |
|---|---|---|---|---|---|---|---|---|---|
| *Ra* | *MRR* | *Ra* | *MRR* | *Ra* | *MRR* | *Ra* | *MRR* | *Ra* | *MRR* |
| 2.23 | 730 | 2.274947 | 485.692 | 2.16589 | 2028.846 | 0.0202 | 0.3347 | 0.0287 | 1.7792 |
| 2.29 | 1460 | 2.317309 | 1534.32 | 2.183 | 1457.984 | 0.0119 | 0.0509 | 0.0467 | 0.0014 |
| 2.32 | 2190 | 2.347448 | 2249.61 | 2.20011 | 887.1224 | 0.0118 | 0.0272 | 0.0517 | 0.5949 |
| 2.37 | 1460 | 2.383297 | 1461.15 | 2.23731 | 1782.861 | 0.0056 | 0.0008 | 0.056 | 0.2211 |
| 2.4 | 2920 | 2.415438 | 2855.25 | 2.26621 | 2915.335 | 0.0064 | 0.0222 | 0.0557 | 0.0016 |
| 2.42 | 4380 | 2.435357 | 3916.01 | 2.29511 | 4047.809 | 0.0063 | 0.1059 | 0.0516 | 0.0758 |
| 2.58 | 2920 | 2.544441 | 3078.74 | 2.38016 | 1290.89 | 0.0138 | 0.0544 | 0.0775 | 0.5579 |
| 2.6 | 5840 | 2.556142 | 5163.77 | 2.43263 | 5830.036 | 0.0169 | 0.1158 | 0.0644 | 0.0017 |
| 2.62 | 5760 | 2.55562 | 6915.46 | 2.4851 | 10369.18 | 0.0246 | 0.2006 | 0.0515 | 0.8002 |
| 1.5 | 1460 | 1.453706 | 1461.6 | 1.5911 | 1785.037 | 0.0309 | 0.0011 | 0.0607 | 0.2226 |
| 1.54 | 2920 | 1.485454 | 3181.96 | 1.6162 | 2925.236 | 0.0354 | 0.0897 | 0.0495 | 0.0018 |
| 1.55 | 4380 | 1.504979 | 4568.98 | 1.64129 | 4065.435 | 0.029 | 0.0431 | 0.0589 | 0.0718 |
| 1.59 | 2920 | 1.546638 | 3119.06 | 1.65388 | 3005.67 | 0.0273 | 0.0682 | 0.0402 | 0.0293 |
| 1.6 | 5840 | 1.579975 | 6072.16 | 1.69076 | 5849.205 | 0.0125 | 0.0398 | 0.0567 | 0.0016 |
| 1.61 | 8760 | 1.601089 | 8691.92 | 1.72764 | 8692.74 | 0.0055 | 0.0078 | 0.0731 | 0.0077 |
| 1.62 | 5840 | 1.676945 | 6100.63 | 1.77943 | 5446.937 | 0.0352 | 0.0446 | 0.0984 | 0.0673 |
| 1.64 | 11680 | 1.713461 | 11519.2 | 1.83988 | 11697.14 | 0.0448 | 0.0138 | 0.1219 | 0.0015 |
| 1.65 | 17520 | 1.737754 | 16604.5 | 1.90033 | 17947.35 | 0.0532 | 0.0523 | 0.1517 | 0.0244 |
| 0.5 | 2920 | 0.505473 | 2781.75 | 0.44881 | 1300.503 | 0.0109 | 0.0473 | 0.1024 | 0.5546 |
| 0.48 | 5840 | 0.516126 | 5837.06 | 0.48976 | 5841.165 | 0.0753 | 0.0005 | 0.0203 | 0.0002 |
| 0.51 | 8760 | 0.514557 | 8559.03 | 0.53072 | 10381.83 | 0.0089 | 0.0229 | 0.0406 | 0.1851 |
| 0.55 | 5840 | 0.567763 | 5794.55 | 0.4944 | 5435.81 | 0.0323 | 0.0078 | 0.1011 | 0.0692 |
| 0.6 | 11680 | 0.603476 | 11845.9 | 0.54714 | 11679.81 | 0.0058 | 0.0142 | 0.0881 | 2E-05 |
| 0.62 | 17520 | 0.626966 | 17564 | 0.59988 | 17923.81 | 0.0112 | 0.0025 | 0.0325 | 0.023 |
| 0.65 | 11680 | 0.636788 | 11486.8 | 0.58558 | 13706.42 | 0.0203 | 0.0165 | 0.0991 | 0.1735 |
| 0.76 | 23360 | 0.72262 | 23530.3 | 0.66189 | 23357.09 | 0.0492 | 0.0073 | 0.1291 | 0.0001 |
| 0.82 | 35040 | 0.79623 | 35240.5 | 0.7382 | 33007.76 | 0.029 | 0.0057 | 0.0998 | 0.058 |

Table 2.3 shows the maximum and minimum error corresponding to *Ra* and *MRR* values for both sets of models. Table 2.4 shows comparison of Mean Absolute Percentage Deviation (MAPD).

**Table 2.3: Comparison of maximum and minimum error**

| Error Type | Models of Eqs. (2.3) and (2.4) | | Models of Eqs. (2.1) and (2.2) | |
|---|---|---|---|---|
| | *Ra* | *MRR* | *Ra* | *MRR* |
| Maximum Error | 2.556142 | 35240.5 | 2.4851 | 33007.76 |
| Minimum Error | 0.505473 | 485.692 | 0.44881 | 887.1224 |

**Table 2.4: Comparison of Mean Absolute Percentage Deviation (MAPD)**

| MAPD for Models of Eqs. (2.3) and (2.4) | | MAPD for Models of Eqs. (2.1) and (2.2) | |
|---|---|---|---|
| *Ra* | *MRR* | *Ra* | *MRR* |
| 0.0235 | 0.0518 | 0.0707 | 0.2047 |

Based on the above information it can be concluded that, Equations (2.3) and (2.4) perform better in representing the system responses. Therefore, these models will be used as optimization criteria.

## 2.3 Definition of design variables

The design variables for face milling operation are the machining parameters, which are:

$Vc$ = Cutting speed, m/min

$fz$ = Feed rate, mm/tooth

$t$ = Depth of cut, mm

## 2.4 Optimization criterion

For the face milling operation, the objective is to select the machining parameters to obtain minimum surface roughness (*Ra*) and maximum Material Removal Rate (*MRR*) simultaneously, through suitable tradeoff between these machining responses. The optimization criteria, i.e. objective functions are the functional relationships among the response variables and the design variables, expressed by Equations (2.3) and (2.4):

$$Ra = 3.082776 - 0.01425 \times Vc + 4.330794 \times fz + 0.465279 \times t + 1.85 \times 10^{-5} \times Vc^2 \\ - 5.78704 \times fz^2 - 0.15278 \times t^2 - 0.00862 \times Vc \times fz - 0.00142 \times Vc \times t - 2.73511 \times fz \times t \\ + 0.018687 \times Vc \times fz \times t$$

$$\begin{aligned}MRR = \ &-2345.09 + 10.37705 \times Vc + 24983.71 \times fz + 7079.912 \times t - 0.01672 \times Vc^2 \\ &- 34722.2 \times fz^2 - 4166.67 \times t^2 - 64.9643 \times Vc \times fz - 13.6425 \times Vc \times t - 66322.5 \times fz \times t \\ &+ 1403.922 \times Vc \times fz \times t\end{aligned}$$

## 2.5 Formulation of constraints

For the optimization problem under study, the restrictions on the design variables (machining parameters) can be formulated into the following constraints:

$$78 - Vc \leq 0, \text{ m/min}; \qquad Vc - 314 \leq 0, \text{ m/min}$$
$$0.04 - fz \leq 0, \text{ mm/tooth}; \qquad fz - 0.16 \leq 0, \text{ mm/tooth}$$
$$0.2 - t \leq 0, \text{ mm}; \qquad t - 0.6 \leq 0, \text{ mm}$$

After formulating the multi-objective optimization problem, five intended optimization routines is going to be employed to solve the problems in Section 3.0.

## 3.0 Multi-objective Optimization Routines

Many different multi-objective optimization routines are available for use in different cases. Multi-objective Optimization problem of face milling operation is going to be solved by using five multi-objective optimization techniques, which are−

i. Global Criterion Method
ii. Lexicographic Method
iii. Weighted Sum Method
iv. ε – Constraint Method
v. Genetic Algorithm

Brief descriptions of these multi-objective optimization techniques are given in Section 3.1 to Section 3.5, mainly derived from Arora, 2012.

## 3.1 Global Criterion Method

Global Criterion Method involves no subjective information on the preference of decision maker regarding the multiple objective functions ($f_i(x)$), where $i = 1$ to $I$. This is a scalarization method that combines all objective functions to form a single function which is then minimized. At first, optimal solution, $f_i(x^*)$ for the individual objective function, $f_i(x)$ is to be determined. The objective function is formulated to minimize the difference between the objective function and the optimal solution for the individual objective function, with respect to the optimal solution, as expressed in Equation (3.1).

Minimize

$$\left\{ \sum_{i=1}^{I} \left\{ \frac{f_i(x^*) - f_i(x)}{f_i(x^*)} \right\}^p \right\}^{\frac{1}{p}} \quad (3.1)$$

subject to

$$g_j(x) \leq 0; \; j = 1 \text{ to } J$$
$$h_k(x) = 0; \; k = 1 \text{ to } K$$
$$lb \leq x \leq ub$$

In the above formulations, $g_j(x)$s are the inequality constraints, $h_k(x)$s are the equality constraints, *lb* represents lower bound and *ub* represents upper bound. The objective function involves a parameter *p*, and the obtained solution depends on its value.

Advantages of Global Criterion Method:

1. It gives a clear interpretation of minimizing the distance from the optimal point, for each individual objective function.
2. The parameter *p* in the objective function can be varied to obtain different optimal solutions and a Pareto front can be generated.
3. For *p* = 2, the objective function is quadratic, which gives a global optimal solution.
4. No subjective information on the preference of decision maker regarding the multiple objective functions is required.

Disadvantages of Global Criterion Method:

1. The use of the optimal point requires minimization of each objective function, which can be computationally expensive.
2. The setting of parameters is not intuitively clear when only one solution point is desired.

## 3.2 Lexicographic Method

In lexicographic method, preferences are imposed by ordering the objective functions according to their attractiveness and importance to the decision maker. Objective functions $f_i(x)$ are ranked such that $f_1(x)$ is most important; then $f_2(x)$, $f_3(x)$ and $f_I(x)$ is of least importance. The optimization problem is formulated as follows:

Minimize

$$f_i(x)$$

subject to

$$f_l(x) \leq f_l(x_l^*); \; l = 1 \text{ to } (i-1); \; i > 1$$
$$g_j(x) \leq 0; \; j = 1 \text{ to } J$$

$$h_k(x) = 0; \quad k = 1 \text{ to } K$$

$$lb \leq x \leq ub$$

This method involves a series of iterations. Here, $i$ represents a function's position in the preferred sequence, and $f_l(x_l^*)$ represents the minimum value for the $l$th objective function, found in the $l$th optimization problem. The algorithm terminates once a unique optimum is determined. Generally, this is indicated when two consecutive optimization problems yield the same solution point.

Advantages of Lexicographic Method:

1. It is a unique approach to specifying preferences.
2. It does not require that the objective functions be normalized.
3. It always provides a Pareto optimal solution.

Disadvantages of Lexicographic Method:

1. It can require the solution of many single-objective problems to obtain just one solution point.
2. It requires that additional constraints to be imposed.

### 3.3 Weighted Sum Method

This is a method that combines all objective functions by multiplying them with certain weighting factor $w_i$ to form a single function which is then minimized. The relative value of the weights generally reflects the relative importance of the objectives. The optimization problem is formulated as follows:

Minimize

$$\sum_{i=1}^{I} w_i f_i(x) \tag{3.2}$$

subject to

$$g_j(x) \leq 0; \quad j = 1 \text{ to } J$$

$$h_k(x) = 0; \quad k = 1 \text{ to } K$$

$$lb \leq x \leq ub$$

$$w_i > 0$$

$$\sum_{i=1}^{I} w_i = 1$$

The weights can be used in two ways:

1. The user may either set $w_i$ to reflect preferences

2. Systematically alter $w_i$ to yield different Pareto optimal points

Advantages of Weighted Sum Method:

1. Easy to use.
2. If all of the weights are positive, the minimum of Eq. (3.2) is always Pareto optimal.
3. Weights can be used either to generate a single solution or multiple solutions.

Disadvantages of Weighted Sum Method:

1. A satisfactory a priori weight selection does not necessarily guarantee that the final solution will be acceptable; one may have to re-solve the problem with new weights. In fact, this is true of most weighted methods.
2. It is impossible to obtain points on non-convex portions of the Pareto optimal set in the criterion space (Marler and Arora, 2010).
3. Varying the weights consistently and continuously may not necessarily result in an even distribution of Pareto optimal points and an accurate, complete representation of the Pareto optimal set.

## 3.4 $\varepsilon$ – Constraint Method

The $\varepsilon$-constraint approach minimizes the single most important objective function $f_s(x)$ with other objective functions treated as constraints. The optimization problem is formulated as follows:

Minimize
$$f_s(x)$$
subject to
$$g_j(x) \leq 0; \; j = 1 \; to \; J$$
$$h_k(x) = 0; \; k = 1 \; to \; K$$
$$lb \leq x \leq ub$$
$$f_i(x) \leq \varepsilon_i; i = 1 \; to \; I; i \neq s$$

Here, $\varepsilon_i$ are the upper bounds for $f_i(x)$. Preferences can be imposed by setting limits on the objectives. Systematic variation of $\varepsilon_i$ yields a set of Pareto optimal solutions. Solutions with active $\varepsilon$-constraints are necessarily Pareto optimal (Carmichael, 1980).

Advantages of ε – Constraint Method:
1. It focuses on a single objective with limits on others.
2. It always provides a weakly Pareto optimal point, assuming that the formulation gives a solution.
3. It is not necessary to normalize the objective functions.
4. It gives Pareto optimal solution if one exists and is unique.

Disadvantages of ε – Constraint Method:
1. The optimization problem may be infeasible if the bounds on the objective functions are not appropriate

## 3.5 Genetic Algorithm

The algorithm begins by creating a random initial population. The algorithm then creates a sequence of new populations. At each step, the algorithm uses the individuals in the current generation to create the next population. To create the new population, the algorithm performs the following steps:

1. Each member of the current population is scored by computing its fitness value.
2. The raw fitness scores are scaled to convert them into a more usable range of values.
3. Selects members, called parents, based on their fitness.
4. Some of the individuals in the current population that have lower fitness are chosen as elite. These elite individuals are passed to the next population.
5. Children are produced from the parents. Children are produced either by making random changes to a single parent—mutation—or by combining the vector entries of a pair of parents—crossover.
6. The current population is replaced with the children to form the next generation.

The algorithm stops when one of the stopping criteria is met. An appeal of genetic algorithm is the ability to converge to the Pareto optimal set rather than a single Pareto optimal point

# 4.0 A Case Study of Multi-objective Optimization: Face Milling

The multi-objective optimization problem that has been formulated in Section 2 is solved using the multi-objective optimization routines discussed in Section 3. Suitable combination of machining parameters – cutting speed ($Vc$), feed rate ($fz$) and depth of cut ($t$) is determined using the optimization routines to obtain optimum machining responses−Material Removal Rate ($MRR$) and surface roughness ($Ra$).

## 4.1 Application of Global Criterion Method

In case of Global Criterion Method, at first optimum values for individual objective functions are determined, which are shown in Table 4.1.

**Table 4.1: Individual optimum values of the objective functions (machining responses)**

|  | Surface roughness, $Ra$ ($\mu$m) | Material Removed Rate, $MRR$ (mm$^3$/min) |
|---|---|---|
| Optimum values | 0.5055 | 35241 |

The obtained optimum values are used to formulate the objective function for Global Criterion Method, minimizing the difference between the objective function and the optimal solution for the individual objective function, with respect to the optimal solution, as expressed in Equation (3.1). Table 4.2 shows set of Pareto optimal solutions which are obtained varying the parameter; $p$. Figure 4.1 shows the Pareto front in criterion space.

**Table 4.2: Set of Pareto optimal solutions obtained in Global Criterion Method**

| Parameter, $p$ | Optimal Solution | | | Optimal Response | |
|---|---|---|---|---|---|
| | Cutting speed, $Vc$ (m/min) | Feed rate, $fz$ (mm/tooth) | Depth of cut, $t$ (mm) | Surface roughness, $Ra$ ($\mu$m) | Material Removed Rate, $MRR$ (mm$^3$/min) |
| 1 | 314 | 0.1050 | 0.6 | 0.6869 | 23047 |
| 2 | 314 | 0.1154 | 0.6 | 0.7111 | 25448 |
| 4 | 314 | 0.1070 | 0.6 | 0.6924 | 23580 |
| 6 | 314 | 0.1060 | 0.6 | 0.6901 | 23358 |
| 8 | 314 | 0.1050 | 0.6 | 0.6878 | 23136 |
| 10 | 314 | 0.1050 | 0.6 | 0.6878 | 23136 |
| 12 | 314 | 0.1040 | 0.6 | 0.6855 | 22914 |
| 14 | 314 | 0.1040 | 0.6 | 0.6855 | 22914 |
| 16 | 314 | 0.1040 | 0.6 | 0.6855 | 22914 |
| 18 | 314 | 0.1040 | 0.6 | 0.6855 | 22914 |
| 20 | 314 | 0.1040 | 0.6 | 0.6855 | 22914 |

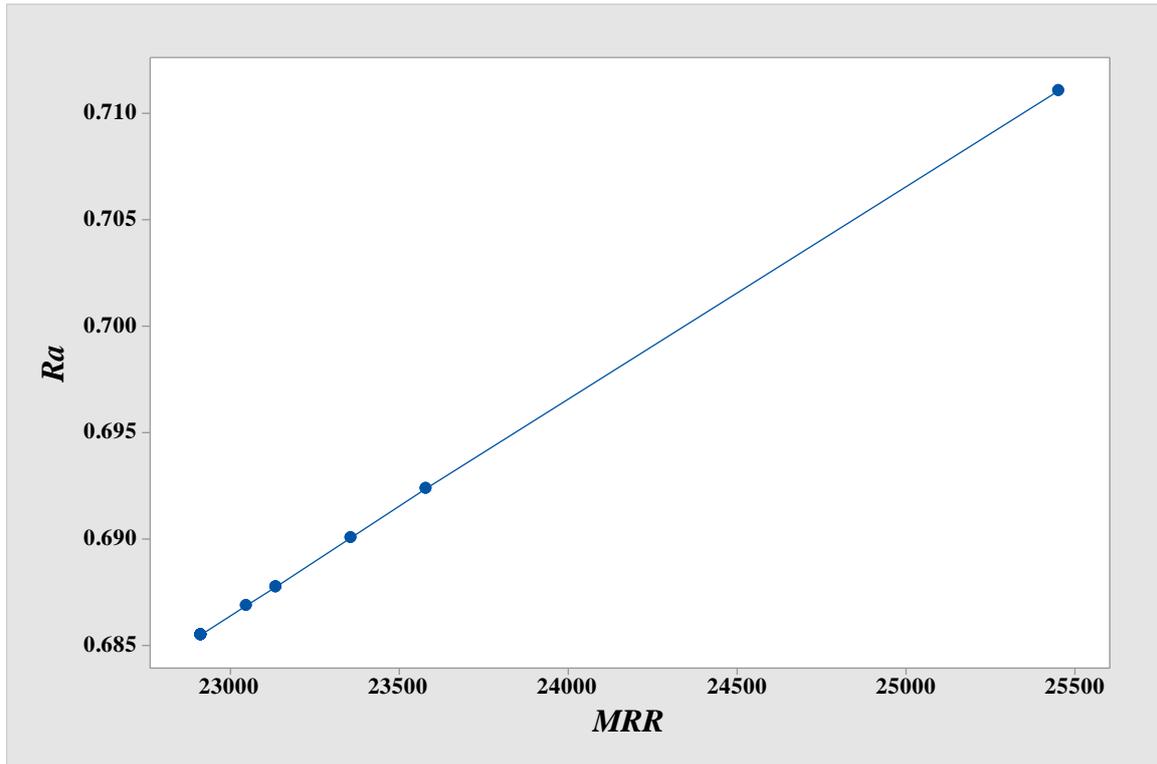

**Figure 4.1: Pareto front obtained in Global Criterion Method**

## 4.2 Application of Lexicographic Method

In case of Lexicographic Method, preferences are imposed by ordering the objective functions according to their attractiveness and importance. In the case study of face milling, importance of *MRR* is set higher compared to *Ra*. Therefore, Equation (2.4) is optimized in the first iteration and optimum value is obtained. Next, Equation (2.3) is optimized in the second iteration with an additional constraint: $MRR \geq MRR_{max}$. Table 4.3 shows the results obtained in two iterations of Lexicographic Method. It is observed that, two iterations give identical solution.

**Table 4.3: Solutions obtained in Lexicographic Method**

| Iteration | Optimal Solution | | | Optimal Response | |
|---|---|---|---|---|---|
| | Cutting speed, *Vc* (m/min) | Feed rate, *fz* (mm/ tooth) | Depth of cut, *t* (mm) | Surface roughness, *Ra* (μm) | Material Removed Rate, *MRR* (mm³/min) |
| 1 | 314 | 0.16 | 0.6 | 0.7962 | 35241 |
| 2 | 314 | 0.16 | 0.6 | 0.7962 | 35241 |

## 4.3 Application of Weighted Sum Method

In case of Weighted Sum Method, the objective function is formulated following Equation (3.2). The value of weighting factor $w_i$ can be systematically altered (e.g. 0 to 1; with increment of 0.1) to yield different Pareto optimal points. One important idea is to normalize the objective functions if they differ highly in magnitude; as in the current case. The objective functions can be normalized using Equation (4.1).

$$f_i^{norm} = \frac{f_i - f_i^{min}}{f_i^{max} - f_i^{min}} \tag{4.1}$$

Table 4.4 shows set of Pareto optimal solutions which are obtained varying $w_i$. Figure 4.2 shows the Pareto front in criterion space.

**Table 4.4: Set of Pareto optimal solutions obtained in Weighted Sum Method**

| Weighting factor, w | Optimal Solution | | | Optimal Response | |
|---|---|---|---|---|---|
| | Cutting speed, $Vc$ (m/min) | Feed rate, $fz$ (mm/tooth) | Depth of cut, $t$ (mm) | Surface roughness, $Ra$ (μm) | Material Removed Rate, $MRR$ (mm³/min) |
| 0 ~ 0.8 | 314 | 0.16 | 0.6 | 0.7962 | 35241 |
| 0.9 | 314 | 0.04 | 0.6 | 0.5149 | 8559 |
| 1.0 | 314 | 0.04 | 0.2 | 0.5055 | 2781.8 |

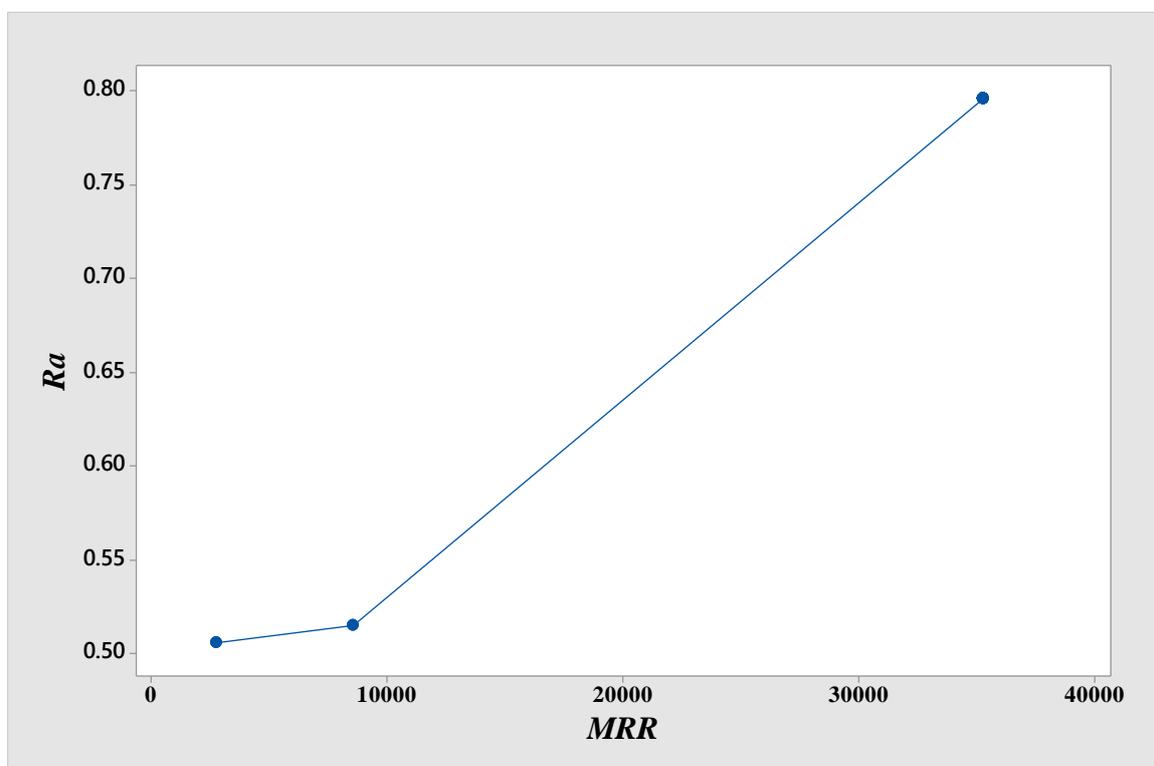

**Figure 4.2: Pareto front obtained in Weighted Sum Method**

## 4.4 Application of ε – Constraint Method

In case of *ε*-Constraints Method, the single most important objective function (Equation (2.4) in current case) is minimized with other objective function (Equation (2.3)) treated as constraints. The upper bound for *Ra* can be altered to yields a set of Pareto optimal solutions. Table 4.5 shows set of Pareto optimal solutions which are obtained varying $\varepsilon_{Ra}$ uniformly between $Ra_{min}$ and $Ra_{max}$ to cover the whole range. Figure 4.3 shows the Pareto front in criterion space.

**Table 4.5: Set of Pareto optimal solutions obtained in ε – Constraint Method**

| $\varepsilon_{Ra}$ | Optimal Solution | | | Optimal Response | |
|---|---|---|---|---|---|
| | Cutting speed, *Vc* (m/min) | Feed rate, *fz* (mm/ tooth) | Depth of cut, *t* (mm) | Surface roughness, *Ra* (μm) | Material Removed Rate, *MRR* (mm³/min) |
| 0.5158 | 314 | 0.0404 | 0.6 | **0.5158** | 8648 |
| 0.7107 | 314 | 0.1153 | 0.6 | **0.7107** | 25409 |
| 0.9159 | 314 | 0.16 | 0.6 | 0.7962 | 35241 |
| 1.1211 | 314 | 0.16 | 0.6 | 0.7962 | 35241 |
| 1.3263 | 314 | 0.16 | 0.6 | 0.7962 | 35241 |
| 1.5315 | 314 | 0.16 | 0.6 | 0.7962 | 35241 |
| 1.7366 | 314 | 0.16 | 0.6 | 0.7962 | 35241 |
| 1.9418 | 314 | 0.16 | 0.6 | 0.7962 | 35241 |
| 2.1470 | 314 | 0.16 | 0.6 | 0.7962 | 35241 |
| 2.3522 | 314 | 0.16 | 0.6 | 0.7962 | 35241 |
| 2.5574 | 314 | 0.16 | 0.6 | 0.7962 | 35241 |

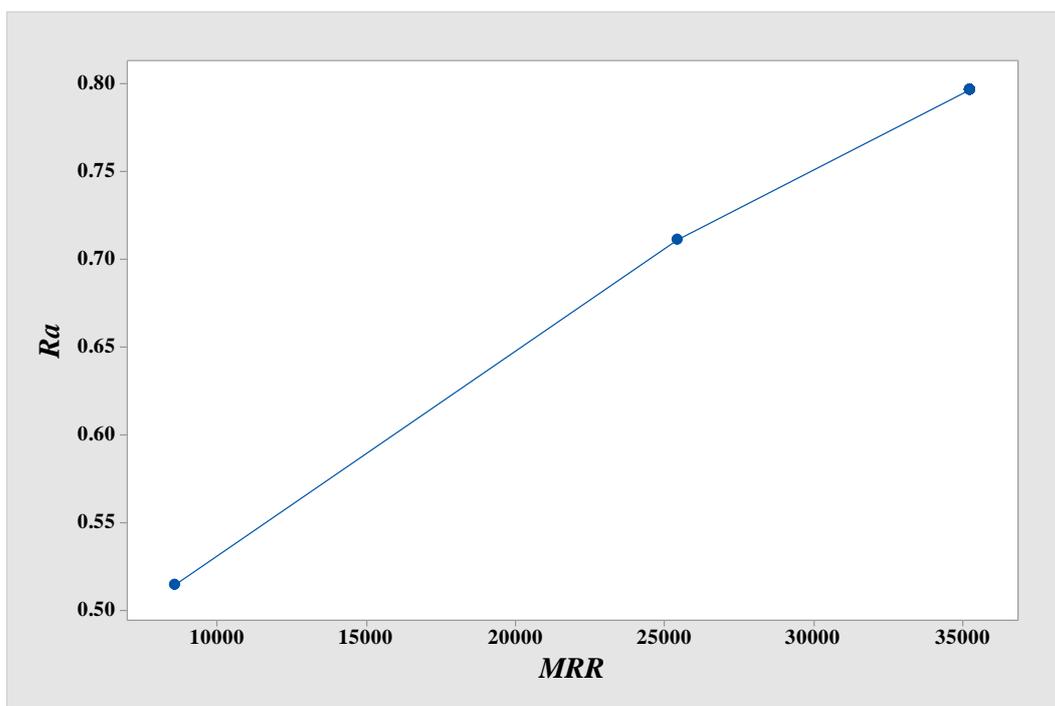

**Figure 4.3: Pareto front obtained in ε – Constraint Method**

## 4.5 Application of Genetic Algorithm

In case of Genetic Algorithm, the multi-objective optimization problem formulated in Section 2 is solved and a set of Pareto optimal solutions is obtained as shown in Table 4.6. Figure 4.4 shows the Pareto front in criterion space.

**Table 4.6: Set of Pareto optimal solutions obtained in Genetic Algorithm**

| Index | Optimal Solution | | | Optimal Response | |
|---|---|---|---|---|---|
| | Cutting speed, $V_c$ (m/min) | Feed rate, $f_z$ (mm/ tooth) | Depth of cut, $t$ (mm) | Surface roughness, $R_a$ (μm) | Material Removed Rate, $MRR$ (mm³/min) |
| 1 | 313.9944 | 0.16 | 0.6 | 0.796247 | 35239.88 |
| 2 | 313.8412 | 0.040241 | 0.574098 | 0.516701 | 8273.971 |
| 3 | 313.3291 | 0.040241 | 0.330112 | 0.516497 | 4829.822 |
| 4 | 313.9944 | 0.16 | 0.6 | 0.796247 | 35239.88 |
| 5 | 313.8223 | 0.056194 | 0.577796 | 0.562547 | 11791.66 |
| 6 | 313.8465 | 0.139075 | 0.591301 | 0.75722 | 30210.68 |
| 7 | 313.4541 | 0.040241 | 0.205112 | 0.508056 | 2879.604 |
| 8 | 313.7007 | 0.083507 | 0.545548 | 0.631147 | 16737.6 |
| 9 | 313.9705 | 0.123172 | 0.584375 | 0.724098 | 26468.48 |
| 10 | 313.904 | 0.093485 | 0.597847 | 0.660724 | 20497.86 |
| 11 | 313.731 | 0.067892 | 0.566768 | 0.593998 | 14062.44 |
| 12 | 313.4881 | 0.041277 | 0.220199 | 0.511083 | 3209.405 |
| 13 | 313.9202 | 0.084704 | 0.570588 | 0.636112 | 17739.54 |
| 14 | 313.8583 | 0.112395 | 0.591803 | 0.703335 | 24437.07 |
| 15 | 313.7041 | 0.062845 | 0.559639 | 0.58043 | 12836.77 |
| 16 | 313.9332 | 0.152551 | 0.559792 | 0.770931 | 31391.3 |
| 17 | 313.8941 | 0.151992 | 0.58674 | 0.778786 | 32753.41 |
| 18 | 313.8412 | 0.102741 | 0.574098 | 0.678937 | 21680.69 |

One limitation of Genetic Algorithm (used in MATLAB) is− it generates set of Pareto optimal solutions which vary from run to run. Therefore, selection of set of Pareto optimal solution is subjective and may vary from user to user. The calculation efficiency of Genetic Algorithm in solving a multi-objective optimization problem cannot be determined with confidence as the number of iterations and number of function count in solving the multi-objective optimization problem vary from run to run.

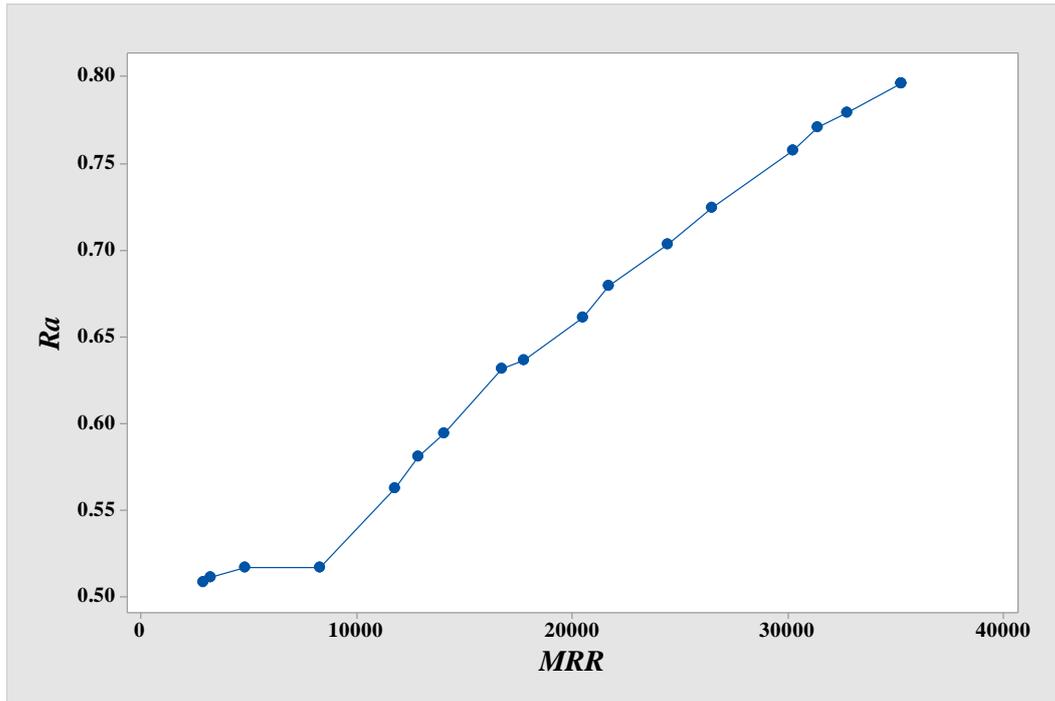

**Figure 4.4: Pareto front obtained in Genetic Algorithm**

## 5.0 Comparison of Multi-objective Optimization Routines

Sections (4.1) to (4.5) discuss the application of five multi-objective optimization routines solving the problem formulated in Section 2. The five routines deal with the multi-objective optimization problem using different approaches. Knowledge of these approaches can help the decision maker use the most appropriate one in solving a multi-objective optimization problem. In current case, among five routines, four (except Lexicographic Method) give set of Pareto optimal solution which may help decision maker choose a suitable tradeoff between levels of *MRR* and *Ra*. Figure 4.5 shows all the Pareto fronts obtained by the four approaches. The computational efficiency of the optimization routines can be measured to find their relative attractiveness to the user. Table 5.1 shows the total number of iterations and function-counts for each multi-objective optimization routines used.

**Table 5.1: Comparison of computational efficiency of the optimization routines**

| Optimization routines | Total Iterations | Total Function counts |
|---|---|---|
| Global Criterion Method | 303 | 1305 |
| Lexicographic Method | 34 | 153 |
| Weighted Sum Method | 165 | 716 |
| ε – Constraint Method | 142 | 620 |
| Genetic Algorithm | 271 | 13601 |

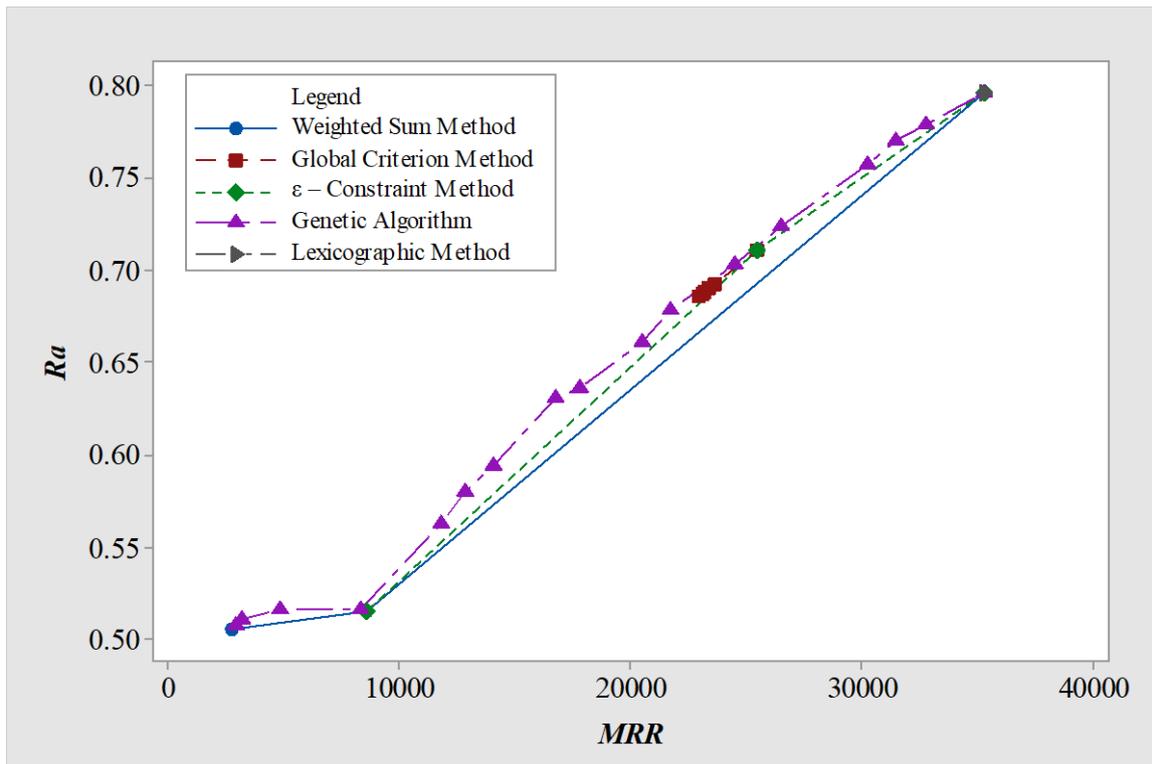

**Figure 4.5: Pareto front obtained in Weighted Sum Method, Global Criterion Method, ε – Constraint Method, Genetic Algorithm and Lexicographic Method**

# 5.0 Conclusion

Different multi-objective optimization routines deal with a multi-objective optimization problem from different perspective. Systematic variation of different parameters in the optimization routines yields a set of Pareto optimal solutions. Applications of five different multi-objective optimization routines have been illustrated with the help of a case study of face milling operation. Face milling operation can be optimized by designing it as a multi-objective optimization problem. Decision maker may choose certain set of machining parameters from the Pareto optimal sets, generated by the five multi-objective optimization routines. From the point of view computational efficiency, relative attractiveness of the optimization routines used in the current case can be ranked in descending order as − Lexicographic Method, ε – Constraint Method, Weighted Sum Method, Global Criterion Method and Genetic Algorithm. Application of these multi-objective optimization routines can be extended to other system optimization involving more than one objective.